\documentclass[preprint,aps,showkeys,12pt]{revtex4}
\usepackage{amsmath,amssymb,euscript,yfonts,psfrag,latexsym,dsfont,graphicx,bbm,color,amstext,wasysym,pdfsync,soul,framed}

\newtheorem{theorem}{Theorem}

\newtheorem{proposition}[theorem]{Proposition}

\newtheorem{remark}[theorem]{Remark}
\newtheorem{problem}[theorem]{Problem}

\def \beq{\begin{equation}}
\def \eeq{\end{equation}}
\def \beqa{\begin{eqnarray}}
\def \eeqa{\end{eqnarray}}
\def \beqan{\begin{eqnarray*}}
\def \eeqan{\end{eqnarray*}}
\def \bea{\begin{eqnarray}}
\def \eea{\end{eqnarray}}

\def \proof {\noindent {\em Proof. }}

\newcommand{\R}{\mathbb{R}}

\newcommand{\N}{\mathbb{N}}
\newcommand{\Z}{\mathbb{Z}}

\newcommand{\Prob}{\mathbb{P}}

\newcommand{\qed}{\hfill $\Box$ \vskip 2ex}

\def \proof {\noindent {\em Proof. }}
\def \qed{\hfill $\Box$ \vskip 2ex}

\begin{document}

\title{A layman's note on a class of frequentist hypothesis testing problems}
\author{Michele Pavon}
\affiliation{Dipartimento di Matematica, Universit\`a di Padova, via
Trieste 63, 35121 Padova, Italy\\ {\tt\small pavon@math.unipd.it}}

%
%
\date{\today}

\begin{abstract} It is observed that for testing between simple hypotheses  where the cost of Type I and Type II errors can be quantified, it is better to let the optimization choose the test size.

\end{abstract}


\keywords{hypothesis testing; Neyman-Pearson; optimization;test size}

\maketitle


\section{Hypothesis testing}
Let $(X,\mathcal F,\mu)$ be a $\sigma$-finite measure space and let $\mathcal P$ be the family of  probability measures $\Prob$ on $(X,\mathcal F)$ which are absolutely continuous with respect to $\mu$ so that, for $A\in\mathcal F$,
$$ \Prob (A)=\int_A p(x)d\mu.
$$
Here $p=d\Prob/d\mu$ is the {\em density} (Radon-Nikodym derivative) of $\Prob$ with respect to $\mu$. We are mostly interested in two cases: The first is when $X$ is a Euclidean space $\R^N$ equipped with the Borel $\sigma$-field and $\mu$ is Lebesgue measure. The second is when $X=\Z^N$ or $X=\N^N$ and $\mu$ is counting measure on all subsets of $X$. This allows us treat probability densities and discrete probability distributions simultaneously. 

Let $\Prob_0, \Prob_1\in\mathcal P$ and let $p_0$ and $p_1$ be the corresponding densities with respect to $\mu$. Let $(X_1,\ldots,X_N)$ be the available sample taking values in $X$. We seek a test $\varphi: X\rightarrow \{0,1\}$ such that, if $(x_1,\ldots,x_N)$ are the observed values, $\varphi(x_1,\ldots,x_N)=0$ if we accept $H_0=\{\Prob_0\}$ and  $\varphi(x_1,\ldots,x_N)=1$ if we accept $H_1=\{\Prob_1\}$. Let $\mathcal C$ be the {\em critical region}, namely the subset of observations $x=(x_1,\ldots,x_N)$ such that $\varphi(x_1,\ldots,x_N)=1$, namely where we reject the null hypothesis, cf. e.g. \cite[Chapter 8]{Ki}.

\section{A class of inference problems}
Consider a simple hypothesis testing problem where we can quantify the cost of each error. Namely, if we reject $H_0$ when it is true we incur the cost $c_0>0$ and if we reject $H_1$ when it is true we incur the cost $c_1>0$. This is the case in many applications such as when, on the basis of a sample, we need to decide whether to halt the production of an item which should meet certain required standards. Both producing a whole stock not meeting the requirements or halting the production process when the requirements are met causes certain quantifiable costs. A type $I$ error occurs with probability $\alpha=\Prob_0(\mathcal C)$ while a type $II$ error occurs  with probability $\beta=\Prob_1(\mathcal C^c)$. It is then natural to try to minimise the cost
$$J (\mathcal C)=c_0\Prob_0(\mathcal C)+c_1\Prob_1(\mathcal C^c).
$$
This is a simple unconstrained optimisation problem which can be formalized as follows.
\begin{problem} \label{1}Find a measurable set $\mathcal C\subset X$ such that the following cost function
$$J(\mathcal C)=c_0\Prob_0(\mathcal C)+c_1\Prob_1(\mathcal C^c)=\int_{\mathcal C}\left[c_0 p_0(x)-c_1p_1(x)\right]d\mu +c_1
$$
is minimised or, equivalently abusing notation, minimize
$$J (\mathds{1}_{\mathcal C})=\int_X\mathds{1}_{\mathcal C}\left[c_0 p_0(x)-c_1p_1(x)\right]d\mu
$$
where $\mathds{1}_{\mathcal C}$ is the indicator function of the set $\mathcal C$. 
\end{problem}

Let us introduce the set 
$$Q=\{f\in L^\infty(X,\mathcal F,\mu) |f:X\rightarrow [0,1]\},
$$
and consider the following ``relaxed" version of Problem \ref{1}:
\begin{problem}\label{2}
$${\rm Minimize}_{f\in Q} J(f),
$$
where
$$\quad J(f)=\int_X f(x)\left[c_0 p_0(x)-c_1p_1(x)\right]d\mu.$$
\end{problem}
Observe that the cost function is {\em linear} in $f$ and $Q$ is convex. Thus, this is a convex optimization problem. We recall a few basic facts from convex optimization. Let $K$ be a convex subset of the vector space $V$, let $F:K\rightarrow \R$ be  convex and let $x_0\in K$. Then, the one-sided directional derivative or hemidifferential of $F$ at $x_0$ in direction $x-x_0$ 
$$F'_+(x_0;x-x_0):=\lim_{\epsilon\searrow 0}\frac{F(x_0+\epsilon(x-x_0))-F(x_0)}{\epsilon}
$$
exists for every $x\in K$ (this is a consequence of the monotonicity of the difference quotients). We record next the characterisation of optimality for convex problems, see e.g. \cite[p.66]{K}.

\begin{theorem}\label{optim}Let $K$ be a convex subset of the vector space $V$ and  let $F:K\rightarrow \R$ be  convex. Then, $x_0\in K$ is a minimum point for $F$ over $K$ if and only if it holds
\begin{equation}\label{suffopt}
F'_+(x_0;x-x_0)\ge 0, \quad \forall x\in K.
\end{equation}
\end{theorem}
We can then apply this result to Problem \ref{2}.
\begin{proposition} The minimum in Problem \ref{1} is attained for 
\begin{equation}\label{optcrit}
\mathcal C^*=\{x\in X |c_0p_0(x)\le c_1p_1(x)\}.
\end{equation}
\end{proposition}
\proof We apply Theorem \ref{optim} to Problem \ref{2} and get that a necessary and sufficient condition for $f^*\in Q$ to be a minimum point of $J(f)$ over $Q$ is
\begin{equation}\label{varineq}
J'(f^*;f-f^*)=\int_X\left[f(x)-f^*(x)\right]\left[c_0 p_0(x)-c_1p_1(x)\right]d\mu \ge 0, \quad \forall f\in Q.
\end{equation}
Observe now that $f^*=\mathds{1}_{\mathcal C^*}$ satisfies (\ref{varineq}). Indeed
\begin{eqnarray}\nonumber \int_X\left[f(x)-\mathds{1}_{\mathcal C^*(x)}\right]\left[c_0 p_0(x)-c_1p_1(x)\right]d\mu\\=\int_{\mathcal C^*}\left[f(x)-\mathds{1}_{\mathcal C^*}(x)\right]\left[c_0 p_0(x)-c_1p_1(x)\right]d\mu+\int_{\mathcal (C^*)^c}\left[f(x)-\mathds{1}_{\mathcal C^*}(x)\right]\left[c_0 p_0(x)-c_1p_1(x)\right]d\mu=\nonumber\\\int_{\mathcal C^*}\left[f(x)-1\right]\left[c_0 p_0(x)-c_1p_1(x)\right]d\mu+\int_{\mathcal (C^*)^c}f(x)\left[c_0 p_0(x)-c_1p_1(x)\right]d\mu\ge 0,
\nonumber
\end{eqnarray}
since both integrals in the last line are nonnegative. Indeed, $f(x)-1\le 0$ and, on $\mathcal C^*$, $c_0 p_0(x)-c_1p_1(x)\le 0$ imply that the integrand in the first integral is nonnegative. The integrand of the second integral is the product of two nonnegative functions and is therefore also nonnegative. Finally, since $f^*=\mathds{1}_{\mathcal C^*}$ is an indicator function, it also solves Problem \ref{1}.
\qed
\begin{remark} We can rewrite the optimal critical region in the familiar form 
\begin{equation}\label{critical}\mathcal C^*=\left\{x\in X |\Lambda(x)\ge \frac{c_0}{c_1}\right\}, \quad \Lambda(x)=\frac{p_1(x)}{p_0(x)}.
\end{equation}
Thus, the ratio of the two costs $c_0/c_1$ plays the role of the multiplier associated to the size constraint in the usual Neyman-Pearson approach. The size of the test and its power, are simply
\begin{equation}\label{errors}\alpha^*=\Prob_0\left(\Lambda(x)\ge \frac{c_0}{c_1}\right), \quad \beta^*=\Prob_1\left(\Lambda(x)\ge \frac{c_0}{c_1}\right).
\end{equation}
\end{remark}
\section{Example}
We illustrate this approach in the simple case of testing the mean of a normal distribution with known variance. Let $\mu$ be Lebesgue measure on $\R$, $p_0=\mathcal N (0,36)$ and $p_1=\mathcal N(1.2,36)$. Suppose $(x_1,\ldots,x_N)$ are the observed values from a random sample  and let $\bar{x}_N=(1/N)\sum_{i=1}^Nx_i$ be the sample mean. Let us fix $\alpha=0.05$ and let $N=100$. Then the optimal Neyman-Pearson test has critical region $\mathcal C _{NP}=\{\bar{x}_{100}\ge 0.987\}$. The corresponding error of the second type is $\beta=0.36$. Since only the ratio $(c_0/c_1)$ matters in the minimisation of Problem \ref{1}, we take from here on $c_1=1$. Thus applying the Neyman-Pearson approach with tests of size $0.05$, we incur the cost
$$J(\mathcal C _{NP})=c_0 (0.05)+0.36.
$$
Next, we compare $J(\mathcal C _{NP})$ with $J(\mathcal C^*)=c_0\alpha^*+\beta^*$, with $\mathcal C^*$ given by (\ref{optcrit}) and $\alpha^*$ and $\beta^*$ given by (\ref{errors}), for different values of $c_0$ and $c_1=1$. We get the results of Table \ref{costs}.

\begin{table}[ht]
\caption{Comparison of costs}
\centering
\begin{tabular}{c|cccc}
\hline\hline
&$J(\mathcal C _{NP})$& & &$J(\mathcal C^*)$\\
\hline
$c_0=1$&$0.05+0.36=0.41$& & &$0.1587+0.1587=0.3174$\\
 $c_0=e$&$2.718\times 0.05+0.36=0.495914$& & &$2.718\times 0.06681+0.30854=0.490129$\\
 $c_0=e^2$&$7.387\times 0.05+0.36=0.7293762$& & & $7.387\times 0.02275+0.5=0.668066171$\\
 $c_0=e^3$&$20.07929\times 0.05+0.36=1.3639645$& & & $20.07929\times 0.00621 + 0.69=0.81469239$\\[1ex]
 \hline
\end{tabular}
\label{costs}
\end{table}
We see that in all cases, as expected since $\mathcal C^*$ gives the minimum cost, fixing $\alpha$ a priori without considering the costs of type I and II errors, leads to a higher cost. The costs are closer when $\alpha^*$ is close to $0.05$. Indeed, if $\alpha^*$ happens to be $0.05$, given the form (\ref{critical}) of $C^*$, we have $C^*=C_{NP}$.

 In conclusion, when the cost of the two errors is known, it appears wiser to let the optimization  determine the size of the test through (\ref{errors}).

\end{document}